# *Constructing and Understanding New and Old Scales on Slide Rules*

*István Szalkai[*], Ace Hoffman[**],*

[*] University of Pannonia, Veszprém, Hungary, szalkai@almos.uni-pannon.hu
[**] The Animated Software Co., USA, rhoffman@animatedsoftware.com

### 2017.06.14., 11:40'

**Abstract**

We discuss the practical problems arising when constructing any (new or old) scales on slide rules ([1],[4],[9]), i.e. realizing the theory [4],[9] in the practice [1]. So, the present article completes the general theory discussed in our previous articles [4],[9]. This might help anyone in planning and realizing (mainly the magnitude and labeling of) new scales on slide rules in the future.

In Sections 1-7 we deal with technical problems, Section 8 is devoted to the relationship among different scales. In the last Section we provide an interesting fact as a surprise to those readers who wish to skip this long article.

**1 What the Scales are ?**

The main idea of any (old or new) slide rule is [4] that we write each value $x \in \mathbb{R}$(real number) on a strip at the (geometric) distance

(0)  $\qquad\qquad\qquad d = u \cdot f(x)$

from the left hand end of the scale $S_1$ where $f$ is a strictly monotonic function (whether increasing or decreasing) and $u$ is a fixed unit, depending on the length of our slide rule and on the extremal values $x_{MAX}$ and $x_{MIN}$ we want to place on the scale.

For example, in the reciprocal scale below, the numbers 10, 5, 2, 0.6 etc. are written at the (geometric) distances $d_{10}=1/10=0.1$, $d_5=1/5=0.2$, $d_2=1/2=0.5$ and $d_{0.6}=1/0.6 \approx 1.67$ (units) from $S_1$. The order of these numbers is reversed, since $f_R(x)=1/x$ is a decreasing function. The geometric unit on this scale is where we have placed the number 1 (more exactly, the distance of the symbols $S_1$ and 1), since $d_1=1/1=1$. One can easily check these distances with an (ordinary) ruler and a pocket calculator. No real number is written at the beginning $S_1$ of this scale, since $d=0=1/x$ holds for no $x$, this fact is denoted by the symbol $\infty$, thinking "$d_\infty=0=1/\infty$".

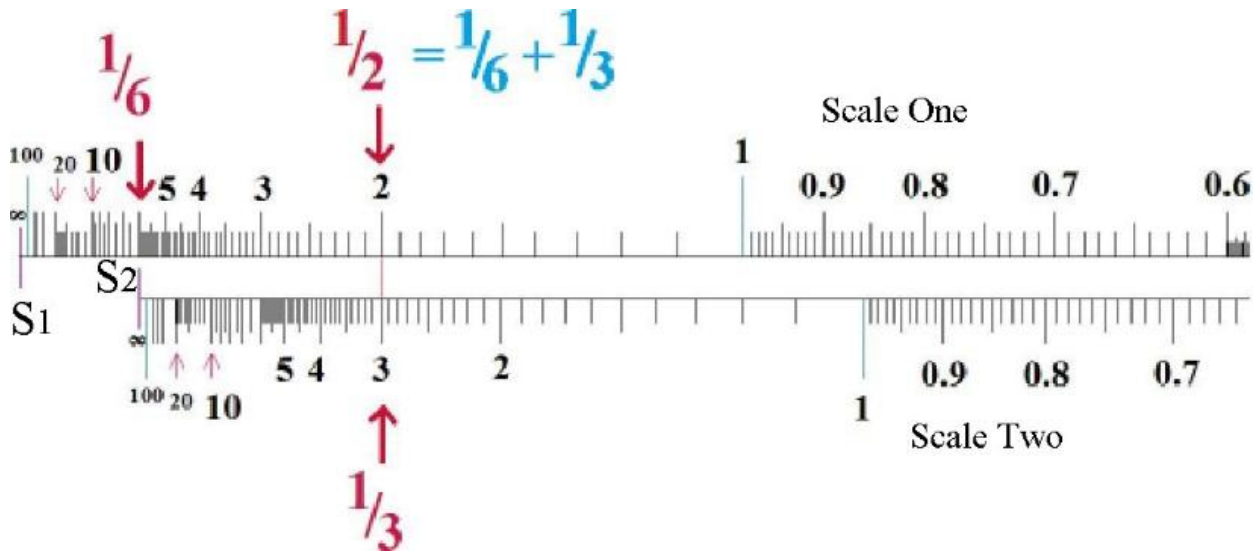

**FIGURE 1. The Reciprocal Scale**

Some other new scales can be seen on Hoffman's downloadable "*Digital Slide Rule*" app and webpage [1] (based on the article [4], the new scales are in green.) The Q scales realize the function $f_Q(x)=x^2$ (quadratic), the $G_1$, $G_2$, $G_4$, $G_5$ scales realize the function $f_G(x) = R \cdot arccos\left(\frac{R}{R+x}\right)$ for $i$=1,2,4,5, where $R$ is the radius of the Earth. This latter function helps us to compute the following everyday problem: if the observer (sailor) at height $h$ (on the mast) observes the object (top of the tower) of height $t$, then their distance is $f_G(h)+f_G(t)$ ([1],[4],[6]). These scales can be easily checked with a ruler and a calculator. The scales $G_3$ and $G_6$ are <u>equidistant</u> scales (that is $f_{G3}(x)=f_{G6}(x)=x$) with the same unit as for scales $G_1$, $G_2$, $G_4$ and $G_5$. This implies that placing the hairline to a number $x$ on scale $G_1$ or $G_4$ (in feet or in meters, respectively), the hairline on scale $G_3$ (in miles) or $G_6$ (in kilometers) shows the length of the horizon we can see from height $x$. The scale on Figure 1 is called "R scale" on Hoffman's *Digital Slide Rule* [1].

The (old, i.e. traditional) base scales C and D of usual slide rules are <u>logarithmic</u> scales, since any number $x$, written on them, has a distance from the left hand end $d_x=log(x)$, i.e., $f_B(x)=log(x)$. This also explains why the C and D scales start with the number $x=1$: its distance is $d_1=log(1)=0$.

## 2 Range of the Scales

Perhaps the first and most important problem is the <u>interval</u> of numbers fit in the scales. Any scale has only some inches length, a small part of the (infinite) real number line. On one scale we can see either the small or the large numbers only, and there is no room for several sizes of one type of any scale. Another problem is when the function $f(x)$ has different slope on different range of the variables $x$ displayed at the same time on the slide rule. For example, let us compare the graph of the function $f(x) = R \cdot arccos\left(\frac{R}{R+x}\right)$ on Figure 2 and the corresponding scale $G_1$ or $G_4$ on Hoffman's *Digital Slide Rule* [1]. For small numbers $f(x)$ is extremely steep and the distances between the numbers are large on the scale, while for larger numbers $f(x)$ is more gentle, where the numbers are closer to each other. This is because each number $x$ is written at distance $d=f(x)$ from $S_1$. (The slope of functions is investigated by derivatives, but this method is cumbersome and in many cases has no clear answer.)

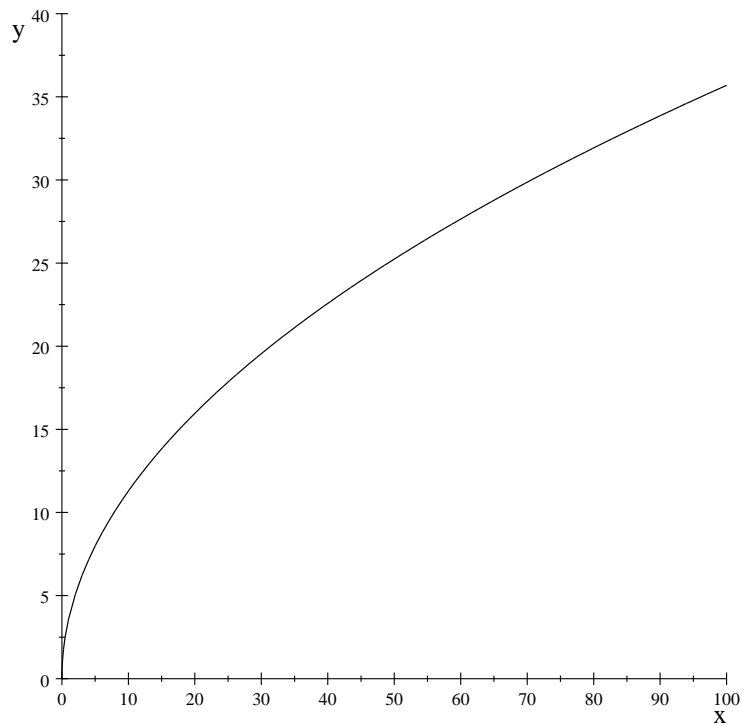

**FIGURE 2.** The Function $R \cdot arccos\left(\frac{R}{R+x}\right)$

The problem of the range of the displayed number was not a problem with the traditional scales, since multiplying any number by $10^p$, i.e. any power of 10, in our calculations, we could fix the slide rule's result to our problem. For example, in multiplication and division or calculating $x^2$, $x^3$ we just have to move the decimal point. Or, calculating logarithm, we just have to add an integer. But, when we calculate Pythagoras' theorem or parallel resistors using our new scales $x^2$ and $1/x$, we have to multiply both $x$ and $y$ with the same number $c$ and this means a multiplication of the result $z$ with the same number $c$. However, this special case does not hold in most of the cases, and never for the scale $R \cdot arccos\left(\frac{R}{R+x}\right)$! Or, what to do when we need a multiple of only $x$ but not of $y$?
Zooming of the scales is possible only on computer animations or on rubber bands. Mathematical calculations for some classes of functions (scales) are discussed in the forthcoming Sections.

### 3 Zooming

Any (new or old) scale has a fixed starting value (in most cases at the left), i.e. $S_1$ has a fixed value, determined by the equality $0 = f(S_1)$. But we may choose the last, right end value $x_{MAX}$ on the scale (despite of its fixed inches length), since the geometric distance on the scales may be zoomed. (Compare small and large slide rules.) Mathematically speaking: for the geometrical distance $d$ we may choose $d=c \cdot f(x)$ for any (positive) real number $c \in \mathbb{R}$. Of course the same zoom factor $c$ must be used for the scales we use jointly in our computations, i.e. one can not use scales $Q_1$ and $Q_2$ in a single movement calculation.

This zooming possibility helps us to choose $x_{MAX}$ freely. For example, on quadratic scales $Q_1$ and $Q_2$ we have $x_{MAX} = 7$ on $Q_1$ and $x_{MAX} = 50$ on $Q_2$. This is essential for functions which have both very steep and almost horizontal parts at the same time on an interval, like $f_G(x) = R \cdot arccos\left(\frac{R}{R+x}\right)$ on Figure 2. For small values $f_G(x)$ is very steep, so the values $x$ on the left hand side on scales $G_1$, $G_2$, $G_4$ and of $G_5$ are rare, but the values are more dense on the right hand side since $f_G(x)$ has a medium slope for larger values. $lim_{x \to \infty} f_G(x) = R$ shows that $f_G(x)$ will be almost horizontal for huge values of $x$, so the values on a scale, containing huge numbers, will be extremely dense. The same phenomenon occurs on traditional scales C and D ($f_C(x)=f_D(x)=log(x)$), called "$x$" scales, or even on scales B and K, called "$x^2$" and "$x^3$" scales ($f_B(x)=(log(x))/2$ and $f_K(x)=(log(x))/3$ since $x_B=(x_C)^2$ and $x_K=(x_C)^3$): the values on the left hand side on these scales are rare, but are more dense on the right hand end. However, this difference of denseness is not so disturbing for us (sometimes remains unobserved), since we get used to writing 10, $10^2$ and $10^3$ etc., that is, our decimal number system resolves this problem.

Finally, we only mention, that scales C, B and K are zoomed images of each other since $f_B(x)=f_C(x)/2$ and $f_K(x)=f_C(x)/3$. Zooming will be an essential tool in the present and our forthcoming articles [5].

### 4 Homogenity

When we shrunk the $Q_1$ (quadratic) scale to go from 0 to 1 and stretched out the $Q_2$ scale (another quadratic) to go from 0 to 10, an interesting thing happened. Namely, that the numbers align themselves, off by one order of magnitude. How can we use two scales (one for the small and the other for large numbers) for our calculations? This experiment suggested to us to investigate the specific property of the scales, called homogenity.

**Definition.** A function $f(x)$ is *homogeneous* if $f(c \cdot x) = c^\alpha \cdot f(x)$ for some exponent $\alpha$ and for all $c \in \mathbb{R}$. □

Clearly all the power functions $f(x) = x^\alpha$ are homogeneous, so all the forthcoming results apply to all the functions $x^\alpha$ ($\alpha=2$ gives the quadratic, $\alpha=-1$ the reciprocal and $\alpha=1$ the *equidistant* scale on your measuring tape). We use $Q_1$ and $Q_2$ to denote the two scales of type $x^\alpha$ we want to compare ($\alpha$ need not be 2 or -1).

**Notations:** **L** is the length of the slide rule, $x_{MAX}^{(1)}$ and $x_{MAX}^{(2)}$ are the maximal numbers printed on the scales $Q_1$ and $Q_2$, **T** $:= x_{MAX}^{(2)}/x_{MAX}^{(1)}$ is the ratio of the maximum numbers, $u_1$, $u_2$ are the units on $Q_1$ and on $Q_2$ respectively, as:

(1) $\qquad d_1(x) = u_1 \cdot x^\alpha \quad \text{and} \quad d_2(x) = u_2 \cdot x^\alpha$

are the geometrical distances where we print $x$ for $0 \leq x \leq x_{MAX}^{(1)}$ on $Q_1$ and for $0 \leq x \leq x_{MAX}^{(2)}$ on $Q_2$. In what follows, we do not write sub- and superscripts when we talk about any of these scales. So

(2) $\qquad L = d_1(x_{MAX}^{(1)}) = u_1 \cdot (x_{MAX}^{(1)})^\alpha = d_2(x_{MAX}^{(2)}) = u_2 \cdot (x_{MAX}^{(2)})^\alpha$

which implies

(3) $\qquad u_1 = L/(x_{MAX}^{(1)})^\alpha \quad \text{and} \quad u_2 = L/(x_{MAX}^{(2)})^\alpha$

(4) $\qquad u_1 \cdot (x_{MAX}^{(1)})^\alpha = u_2 \cdot (x_{MAX}^{(2)})^\alpha = u_2 \cdot (T \cdot x_{MAX}^{(1)})^\alpha$

(5) $\qquad u_1 = u_2 \cdot T^\alpha$.

Now, let us calculate in what conditions are the numbers $x_1$ on $Q_1$ and $x_2$ on $Q_2$ printed on the same place, i.e. on the same distance from $S_1$ :

$$d_1(x_1) = d_2(x_2) \leftrightarrow u_1 \cdot (x_1)^\alpha = u_2 \cdot (x_2)^\alpha$$

$$\leftrightarrow u_2 \cdot T^\alpha \cdot (x_1)^\alpha = u_2 \cdot (T \cdot x_1)^\alpha = u_2 \cdot (x_2)^\alpha \leftrightarrow T \cdot x_1 = x_2 \leftrightarrow T = x_2/x_1 .$$

**Summarizing:** the ratio

(6) $$T := x_{MAX}^{(2)} / x_{MAX}^{(1)}$$

determines which (smaller) numbers $x_2$ and $x_1$ are on the same place on $Q_2$ and on $Q_1$, namely if and only if

(7) $$x_2/x_1 = T = x_{MAX}^{(2)} / x_{MAX}^{(1)} .$$

This phenomenon can be easily observed when T is a power of 10 because of our decimal number system. In other words: if T=10, T=2, T=5 or an "easy" rational number, then there is no significant difference between the scales $Q_1$ and $Q_2$, any of them does the others' job. We say this since we can find out (and use) <u>any</u> of the numbers $x$ on $Q_1$ and $y$ on $Q_2$, $y$ is aligned to(directly above) $x$, since, e.g., $y=2 \cdot x$ in the case T=2. Since all the interstitial nibs are on the same place on both scales, the two scales are <u>equivalent</u>.

So, we are forced to choose a "more irrational" ratio T. One can try values like T=1.3, 1.7, 77.46 ($\approx 10\sqrt{60}$), 83.67 ($\approx 10\sqrt{70}$), … . Since T depends on $\alpha$ and on the tasks we want to solve by our scales, it is impossible to give an easy and general answer. Instead, we have to try out different values of T on the screen (see [1],[7]). Of course these problems are valid for all scales $x^\alpha$ (square, reciprocal, equidistant,etc.). This is why the scales $R_i$ and $Q_i$ have such ugly ranges in [1]. (Let us mention the nice curiosity and pedagogical advantage of scale R that contains values <u>less and greater</u> than 1 at the same time on [1]!)

### 5 Accuracy

The number of correct decimal digits is limited by the length of the scales and the (traditional) materials of the slide rules, and mainly by our eyes. This problem is in contradiction to Section 2 "Range of the Scales", the balance could be approached with the method of Section 3 "Zooming". Here we try to help with some math.

The criterion for reading <u>two</u> decimal values is that the marks $x$ and $1.01 \cdot x$ could be distinguished, i.e. their distance should be at least some given $h$ (depending on our eyes):

(8) $$d(1.01 \cdot x) - d(x) = u \cdot ( f(1.01 \cdot x) - f(x) ) \geq h.$$

For <u>homogeneous</u> functions (8) can be written easier as

(9) $$u \cdot (1.01^\alpha - 1) \cdot f(x) \geq h .$$

Now, if the function $f(x)$ is strictly monotone increasing OR decreasing, we have to check only

(10) $$u \cdot (1.01^\alpha - 1) \cdot f(x_{MIN}) \geq h \quad \text{OR} \quad u \cdot (1.01^\alpha - 1) \cdot f(x_{MAX}) \geq h ,$$

i.e.

(11) $$u \geq \frac{h}{(1.01^\alpha - 1) \cdot f(x_{MIN})} \quad \text{OR} \quad u \geq \frac{h}{(1.01^\alpha - 1) \cdot f(x_{MAX})} ,$$

or, in other words

(12) $$x_{MIN} \geq f^{-1}\left(\frac{h}{(1.01^\alpha - 1) \cdot u}\right) \quad \text{OR} \quad x_{MAX} \leq f^{-1}\left(\frac{h}{(1.01^\alpha - 1) \cdot u}\right) .$$

For the functions $f(x) = x^\alpha$ (11) and (12) mean

(13) $$u \geq \frac{h}{(1.01^\alpha - 1) \cdot (x_{MIN})^\alpha} \text{ if } \alpha > 0 \quad \text{AND} \quad u \geq \frac{h}{(1.01^\alpha - 1) \cdot (x_{MAX})^\alpha} \text{ if } \alpha < 0 ,$$

and

(14) $$x_{MIN} \geq \sqrt[\alpha]{\frac{h}{(1.01^\alpha - 1) \cdot u}} \text{ if } \alpha > 0 \quad \text{AND} \quad x_{MAX} \leq \sqrt[\alpha]{\frac{h}{(1.01^\alpha - 1) \cdot u}} \text{ if } \alpha < 0$$

where $x_{MIN}$ and $x_{MAX}$ are the smallest and largest values we are able to distinguish from its neighbor on the scale.
**Note:** Looking on the scales ([2]) we can observe that the distances between consecutive marks keeps increasing on Q (since $\alpha = 2 > 0$) and decreasing on R ($\alpha = -1 < 0$), according to the above calculations.

## 6 Triangles and the $x^2$ scale

The quadratic scale ($\alpha=2$) has distinguished interest of ours, probably since Pythagoras' theorem and triangles are commonly used and triangles are easy to visualize. That is why we show a short closer look at scale Q, i.e. $f(x)=x^2$. What kind of scales do we need, concerning length and range vs. accuracy for the result of the calculation

(15) $$c = \sqrt{a^2 + b^2},$$

We have to ensure that all (or rather, most of) the calculations in (15) could be done on <u>one</u> scale, i.e. all the three numbers a,b and c must be found on this scale. By the homogenity we are interested in positive numbers $\tau_1<\tau_2$ such that, assuming

(16) $$\tau_1 \cdot a \leq b \leq \tau_2 \cdot a$$
   i.e.
(17) $$\tau_1 \leq b/a \leq \tau_2$$

the equality (15) could be calculated with a single quadratic scale with desired accuracy. (17) is equivalent to

(18) $$\tau_1 \leq \tan(\alpha) \leq \tau_2 \quad \text{i.e.} \quad \arctan(\tau_1) \leq \alpha \leq \arctan(\tau_2)$$

where $\alpha$ is the angle conjugate to the side a. The choice of $\tau_1$ and $\tau_2$ will be clearer later. (17) and (15) imply

(19) $$a \cdot \sqrt{1 + (\tau_1)^2} \leq c \leq a \cdot \sqrt{1 + (\tau_2)^2}.$$

By (15) through (19) our task is to find suitable numbers a, $\tau_1$ and $\tau_2$ so that all the numbers in the intervals

(20) $$[\tau_1 \cdot a, \tau_2 \cdot a] \quad \text{and} \quad [a \cdot \sqrt{1 + (\tau_1)^2}, a \cdot \sqrt{1 + (\tau_2)^2}]$$

could be seen with desired accuracy in a single scale on the Slide Rule. Then (20) ensures that (15) would be solvable by <u>this</u> scale for all triangles, satisfying(17), or equivalently(18). Clearly we also have

(21) $$\tau_1 < \sqrt{1 + (\tau_1)^2} \quad \text{and} \quad \tau_2 < \sqrt{1 + (\tau_2)^2}.$$

## 7 Summary and Examples

 **(i)** (11) determines the smallest/largest numbers that are visible on our scale. Other numbers' accuracy is ensured by the monotonicity of the function $f(x)$.
 **(ii)** We have to choose the numbers a, $\tau_1$ and $\tau_2$, these values determine the intervals in (20) in which (15) is solvable by this scale for all triangles, satisfying (17), or equivalently (18).
**(iii)** Do not choose T in (6) to a power of 10 or any "easy rational number" to avoid the other scale $Q_2$ to be useless.

**Examples:** Let L:=250 mm, h:=0.5 mm and $x_{MAX}$:=100. By (3) u=L/$(x_{MAX})^2$ =250/100² =0.025 and by (11) we have $x_{MIN}=\sqrt{(0.5/(0.025 \cdot 0.0201))} \approx 31.54$. This is the range of our scale: [31.54, 100].
What triangles can be solved on this scale? Let us examine e.g. a=40. If we substitute the interval [$x_{MIN}$, $x_{MAX}$] into (20), we get

(22) $$40 \cdot \tau_1 = 31.54 \implies \tau_1 = 31.54/40 \approx 0.7886 \implies \alpha_1 = \arctan(\tau_1) \approx 38.26°,$$
and
(23) $$40 \cdot \sqrt{1 + (\tau_2)^2} = 100 \implies \tau_2 = \sqrt{\left(\frac{100}{40}\right)^2 - 1} \approx 2.291 \implies \alpha_2 = \arctan(\tau_2) \approx 66.42°.$$

Therefore, if one of the perpendicular sides of the triangle is length a=40, then the other can be

(24) $$a \cdot \tau_1 \approx 31.54 \leq b \leq a \cdot \tau_1 \approx 91.64$$

and the angle $\alpha$ opposite to a, by (18), is

(25) $$38.26° \leq \alpha \leq 66.42°.$$

If we choose a=32 (the smallest, extreme value) on this scale, then we get

(26) $$32 \cdot \sqrt{1 + (\tau_2)^2} = 100 \implies \tau_2 \approx 2.961 \implies \alpha_2 \approx 71.34°.$$

## 8 Naming the Scales

The standard use of the letters A,B,...,Z is the best solution to uniquely distingish and refer to the different scales of slide rules [2],[3]. (We are happy that the letters G, Q and R are generally left empty, so we could use them for our new scales.)

Although the formulas which usually are printed on the other end of the scales, like $e^x$ or $\sqrt{1-x^2}$, are shown to describe the real structure and the use of the scales, we have to rethink them a bit, for understanding the old, "**traditional**" [2],[3] and new scales.

Scales C and D [2],[3] are called "$x$ -scale", B is "$x^2$ -scale", L is "$log(x)$ -scale", LL3 is "$e^x$ -scale", etc. This means, that placing the hairline to a number $x$ in C or D, the same hairline immediately shows $x^2$ on B, $log(x)$ on L, $e^x$ on LL3 and so on. Do not mix this (useful and real) suggestion to the fact, that on these scales the distance-functions are, in fact $f_C(x)=f_D(x)=log(x)$, $f_B(x)=log(\sqrt{x})=log(x)/2$, $f_{LL3}(x)=log(log(x))$ and $f_L(x)=x$. On the basis of the distance functions these scales could perhaps more correctly be called "$log(x)$", "$log(x)/2$", "$log(log(x))$" and "$x$" - scales, but clearly this new ontology would be very confusing. Look on your slide rule: scale L is an equidistant one like an ordinary ruler (measuring tape), this fact is explained by $f_L(x)=x$.

However, the **new** scales (**Q,R,$G_i$ in green** [1]) do not fit into this "traditional" system. It is important to recognize that the values $x$ on our new scales has no connection to the old scales C or D at all. (The only exception among new and old scales is introduced in the last Section.) Summarizing: at this moment old (traditional) and new scales ([1],[4]) should be used separately. Moreover, this is true about the connection among new scales themselves, too. This raised the (not accepted) idea of using different variable letter for the new scales, but then there would be more scales than letters in the alphabet, and all letters are the same for scientists and engineers. This is why we used different colours in [1] and called these scales on their real names as "$x^2$", "$1/x$" and "$R \cdot arccos\left(\frac{R}{R+x}\right)$" scales. (No matter that $Q_1$, $Q_2$, $R_1$,$R_2$ and $G_1$,$G_2$ use different units (zooming).) The interested Reader may observe in Hoffman's *Digital Slide Rule* [1] that only the scales $G_1$, $G_2$, $G_4$ and $G_5$ are "$R \cdot arccos\left(\frac{R}{R+x}\right)$" scales while the scales $G_3$ and $G_6$ are, in fact equidistant scales as the L scale is. (The use of the $G_i$ family of scales is described in the pop up help screen in [1].) This latter fact would suggest to identify the scales $G_3$, $G_6$ and L, but the obligatory different ranges of these three scales deny this possibility. Regarding the L scale saying "$log(x)$" we didn't make that up, it seems to be pretty standard, as are the other equations which actually use logs but don't say it. Confusing for the novice, of course, but at least it matches standard industry practice!

Finally let us highlight, that new scales can solve some calculations (the horizont problem, Pythagoras, parallel resistors, etc. [1],[4]) much faster, and in just one movement, than without them.

## 9 Miscellaneous

Many other technical problems arose which we do not have space here to discuss. The articles [5], [9] contain more mathematical calculations. Finding the most appropriate ranges and zooming of the scales needed much experimention and patience. If one scale is not enough for the desired range (e.g. $G_i$, $Q_i$, $R_i$), how many scales do we need and how to choose their range? What about extendable rubbers? Where to put the many new scales (two or more sided slide rule do we (or a carpenter/metal worker) construct)? On computers we can easily imagine buttons to change the scales, in "real life" such a task would be hard (or impossible) to realize ([1],[7]). Is it possible to put infinite scales on a slide rule [5] ? What to write to the short help screen or info on a small tag or button next to the slide rule [1]?

## 10 An Interesting Phenomenon

We emphasized in Section 8 that there is no connection among new and traditional scales. Here we introduce the only exception we have found so far, however further connections might turn to daylight sometime.

Let us consider the special case when the number $x_R=1$ on the right hand end of the R scale is placed exactly above $x_C=10$ on the C scale. Now, comparing the geometrical distances from the common $S_1$ left to any fixed (geometrical) point we have (using the same unit $u$ for $f_C$ and $f_R$ in (0)):

(26) $\qquad f_C(x_C) = log_{10}(x_C) = 1/x_R = f_R(x_R)$ ,

i.e.

(27) $\qquad x_C = 10^{(1/x_R)} = exp_{10}(1/x_R)$

or

(28) $\qquad x_R = 1/log_{10}(x_C) = log_{x_C}(10)$

where $x_C$ and $x_R$ denote the numbers we see at the same hairline on scales C and R, respectively and $log_{x_C}$ denotes the logarithm of base $x_C$ ($x_C \neq 1$). Functional connections like in (28) are rare but do exist in practical use, namely in computer science, see e.g. Wagstaff-Norman-Campbell's paper [8]. Let us write here some examples and ask the Reader to check them by a slide rule [1].

| $x_C$ | 1.496 | 4.000 | 4.976 | 10.000 , |
|---|---|---|---|---|
| $x_R$ | 5.714 | 1.661 | 1.436 | 1.000 . |

(Recall: $x_C$ should be found on scale C or D, while $x_R$ is the value written on scale R.)